\documentclass{amsart}
\usepackage{amsmath,amssymb,amscd,amsthm}
\usepackage[all]{xy}

\newcommand{\N}         {\mathbb N}
\newcommand{\Z}         {\mathbb Z}
\newcommand{\Q}         {\mathbb Q}
\newcommand{\R}         {\mathbb R}
\newcommand{\C}         {\mathbb C}
\newcommand{\F}         {\mathcal{F}}

\newcommand{\CP}        {\mathbb{CP}}
\newcommand{\Cone}      {\mathrm{Cone}}

\newcommand{\Pf}        {\mathrm{Pf}}

\newcommand{\ra}        {\rightarrow}

\newcommand{\tensor}    {\otimes}

\newcommand{\Union}{\bigcup}
\newcommand{\intersect}{\cap}
\newcommand{\iso}       {\cong}

\DeclareMathOperator{\Sign}{Sign}

\DeclareMathOperator{\Ann}{Ann}
\DeclareMathOperator{\Aff}{Aff}
\DeclareMathOperator{\Lin}{Lin}
\DeclareMathOperator{\Reg}{Reg}

\newcommand{\h}         {\mathfrak{h}}

\newcommand{\tlie}      {\mathfrak{t}}

\newcommand{\hs}        {\mathfrak{h}^*}

\newcommand{\ts}        {\mathfrak{t}^*}

\newtheorem{thm}{Theorem}
\newtheorem{theorem}{Theorem}
\newtheorem{Def}{Definition}
\newtheorem{prop}{Proposition}
\newtheorem{cor}{Corollary}
\newtheorem{lemma}{Lemma}

\begin{document}

\title{A Wall-crossing Formula for the Signature of Symplectic Quotients}

\author{David S. Metzler}

\maketitle

\begin{abstract}
We use symplectic cobordism, and the localization result of 
Ginz\-burg, Guil\-lemin, and Kar\-shon \cite{GGK:1996C} to find a 
wall-crossing formula for the signature of regular symplectic quotients
of Hamiltonian torus actions. The formula is recursive, depending 
ultimately on fixed point data. In the case of a circle action, we obtain
a formula for the signature of singular quotients as well.
We also show how formulas for
the Poincar\'e polynomial and the Euler characteristic (equivalent to those
of Kirwan \cite{FK:1984CQ}) can be expressed
in the same recursive manner.
\end{abstract}

\section{Introduction}
\label{introchap}

Since the original convexity theorems of Atiyah \cite{A:1982CCH} and 
Guillemin-Sternberg \cite{GS:1982C}
on the image of the moment map associated to a Hamiltonian torus action,
it has become clear that the moment polytope is a rich source
of combinatorial invariants of these actions (see \cite{G:1994MM}). 
In particular, traditional topological invariants
of the various reduced spaces associated to the action turn
out to have formulas that involve fixed point data combined with
combinatorial data from the moment polytope. We will simply use the
word ``combinatorial'' to refer to formulas of this type. 
In this paper we focus primarily on the signature, and find formulas for the
signature of reduced spaces. In the circle case (Section 
\ref{sec:circle}) we are
able to find formulas both for regular and singular reductions,
thanks to recent work of Lerman and Tolman. In the torus case
(Section \ref{recursivechap})
we find a recursive formula for the signatures of various
``subreductions,'' including all regular reductions.

Our fundamental object of study is a Hamiltonian $T$-space $M$,
with symplectic form $\omega$ and moment map $\phi: M \ra \ts$.
We denote such a Hamiltonian $T$-space by $(M,\omega,\phi)$.
We allow our spaces to be orbifolds as well as manifolds.
We will usually restrict ourselves to compact spaces. 

Given a point $a \in \ts$, the symplectic reduction, or
symplectic quotient of $M$ at $a$, is defined to be 
  \[
    M_{a} = \phi^{-1}(a)/T,
  \]
sometimes denoted $(M/ \! \! / T)_{a}$ (when we need to explicitly
include $T$).
Recall that if $a$ is a regular value of $\phi$, $M_{a}$ is a smooth
symplectic orbifold, and that the symplectic form $\omega_{a}$ on $M_{a}$ 
is defined by the requirement that 
\[
  \pi^{*}\omega_{a} = i^{*}\omega
\]
where $\pi: \phi^{-1}(a) \ra M_{a}$ and $i: \phi^{-1}(a) \ra M$ are
the projection and inclusion maps, respectively. 

Let $\Delta_{\mathrm{reg}}$ be the set of regular values of the
moment map. We will see in section \ref{xraychap} that this is a finite
union of open convex polytopes, which we call \textit{chambers}.
When computing topological invariants of reduced spaces, the essential fact is
\begin{prop}\label{invarianceofreduction}
  Let $(M,\omega,\phi)$ be a compact, connected Hamiltonian $T$-space. Let
  $a,b \in \Delta_{\mathrm{reg}}$ lie in the same chamber $P$.
  Then the two reduced spaces $M_{a}$ and $M_{b}$ are 
  diffeomorphic. 
\end{prop}
\begin{proof}
  The restricted map $\phi: \phi^{-1}(P) \ra P$ is a submersion
  with compact fibers. By Ehresmann's theorem, 
  it is a fibration. Since $P$ is connected, all of the fibers are
  diffeomorphic. Since $\phi$ is $T$-equivariant, the fibers
  $\phi^{-1}(a)$ and $\phi^{-1}(b)$ are equivariantly diffeomorphic,
  and therefore 
  \[
    M_{a} = \phi^{-1}(a)/T \iso \phi^{-1}(b)/T = M_{b}.
  \]
\end{proof}
Hence any topological invariant of the reduced spaces
$M_{a}$ ($a \in \Delta_{reg}$) will be a function only of
the chamber that $a$ is in, making it
essentially combinatorial.

In section \ref{sec:circle} we use the cobordism results of 
Guillemin, Ginzburg and Karshon \cite{GGK:1996C} to derive a 
formula for the signature of a regular reduced space in the case of
a circle action, in terms of fixed point data. This formula
is best expressed as a ``wall-crossing'' formula. This form
is quite simple, and generalizes best to the torus case.
In Section \ref{subsec:singularreduction} we
analyze the case of singular reduction by a circle action,
using recent results of Lerman and Tolman 
\cite{LermanTolman:1998Sing}.
Again we obtain a simple formula in terms of fixed point data.

Let us make precise what we mean by ``fixed point data.''
If $M$ is a compact $T$-space, the fixed point set $M^{T}$
has finitely many connected components $F_{1},\dots ,F_{m}$.
Each $F_{r}$ is a compact connected manifold. If $M$ is 
Hamiltonian, it has associated to it a natural isotopy class of
$T$-invariant almost complex structures \cite{McDuffSalamon:1995}. 
Hence the normal bundle $N_{r}$ to each $F_{r}$ has a natural
isotopy class of complex structures, which allows us to 
define the weights of the action of $T$ on 
the normal bundle. Call these weights $\alpha_{r,k} \in \ts$,
$k= 1,\dots,\mathrm{rank} (N_{r})/2$. It is easy
to see that the moment map $\phi$ must be constant on each $F_{r}$, so
each component $F_{r}$ defines a point in $\ts$, $\phi_{r} = \phi(F_{r})$.
The manifolds $F_{r}$,
together with the weights $\alpha_{r,k}$ and the points $\phi_{r}$,
comprise the
fixed point data that our formulas ultimately depend on.

To investigate the torus case, 
we first introduce a variant of the \textit{X-ray} of Tolman 
\cite{ST:1996NK}, and describe it in Section \ref{xraychap}. 
This is a package of data that includes the fixed point data,
but also takes into account the fixed point sets of subtori $H \subset T$.
Using this language, in Section \ref{recursivechap}, 
we find a recursive generalization
of the wall-crossing formula obtained in the circle case, where
the data involved in crossing a wall of dimension $d$ is obtained
by a wall-crossing procedure in dimension $d-1$.
In some cases this reduces to a formula which is as simple
as the one in the circle case.

We conclude by deriving similar recursive formulas for the Poincar\'e
polynomial (and hence the Euler characteristic). These formulas
are perhaps more gemoetrical and easier to calculate than 
those of Kirwan \cite{FK:1984CQ}.

This paper is largely derived from the author's Ph. D. thesis
\cite{Metzler:1997}, supervised by Victor Guillemin. The author
would also like to gratefully acknowledge the support and
advice of Eugene Lerman, Sue Tolman, and Yael Karshon.


\section{Circle Actions}
\label{sec:circle}

Since the signature is a cobordism invariant, it is natural to 
use cobordism arguments to try to calculate it.
Guillemin, Ginzburg, and Karshon \cite{GGK:1996C} 
have introduced a notion of
cobordism of spaces with Hamiltonian group actions which is perfectly
suited to our task. It has two main features: first, a cobordism 
of Hamiltonian spaces induces cobordisms
of their reduced spaces. Second, there is a localization theorem, which says
that every space is cobordant to a sum of local models, determined by the
fixed point data. Further, 
it was noticed by Karshon \cite{YK:1997MM} that the ``moment map''
that appears in their definition of cobordism need not actually
come from a symplectic form, but must only satisfy some simple
axioms. This will prove useful in investigating the signature of
singular reductions.

Note: as has become standard in studying symplectic reduction, we
will deal everywhere with orbifolds, since they arise naturally
from reduction anyway. For our purposes this introduces essentially
no difficulties.

\subsection{Abstract Moment Maps and Cobordism}
\label{subsec:abstractmoment}

Let $M$ be an orbifold with an action of a torus $T$. 
Given a subgroup $H \subset T$, let $\pi_{H}: \ts \ra \hs$
denote the natural projection map, and let $M^{H}$ be the
set of fixed points of $H$.
\begin{Def}\label{def:abstractmoment}
  An \textbf{abstract moment map} for $M$ is a $T$-invariant
  map $\phi: M \ra \ts$ such that, for every subgroup $H \subset T$,
  the composed map 
  \[
    \pi_{H} \circ \phi|_{M^{H}}: M^{H} \ra \ts \ra \hs
  \]
  is locally constant. 
\end{Def}

It is easy to show that if $(M,\omega,\phi)$ is a Hamiltonian $T$-space,
then the moment map $\phi$ is an abstract moment map. However,
abstract moment maps are much looser objects (for example, the zero map
is trivially an abstract moment map) and they do not depend on a
symplectic form for their definition. However, the notion of 
reduction still makes sense: it is easy to see that if $a \in \ts$
is a regular value of $\phi$, then the action of $T$ on $\phi^{-1}(a)$
is locally free, and hence $\phi^{-1}(a)/T$ is an orbifold, denoted
by $M_{a}$.

We will need to be careful about orientations of reduced spaces.
In the symplectic case, the reduced space is clearly oriented,
being symplectic. However, even in the case of an abstract moment map,
the reduced space has a canonical orientation, as long as $M$
is itself oriented. For, we have short exact sequences
\[
  \xymatrix@=6pt{ 0 \ar[r] & {\tlie} \ar[r]^{\iota} 
                    & {T_{p} (\phi^{-1}(a))} \ar[r] 
                    & {T_{[p]} (M_{a})} \ar[r] & 0 \\
                  0 \ar[r] & {T_{p} (\phi^{-1}(a))} \ar[r] 
                    & {T_{p}} M \ar[r]^{d \phi} 
                    & {\ts} \ar[r] & 0. }
\]
where $\iota$ denotes the infinitesimal action of the torus.
Hence if we choose an orientation for $\tlie$ (and hence for $\ts$)
we obtain orientations on $\phi^{-1}(a)$ and $M_{a}$. If we
reverse the orientation of $\tlie$, the orientation on
$\ts$ also reverses, and hence the orientation on $M_{a}$
is unchanged, so it is actually well-defined independent of the
orientation of $\tlie$.

We can now introduce the notion of cobordism that we need.
Let $(M_{1},\phi_{1})$ and $(M_{2},\phi_{2})$ be two oriented $T$-spaces
with \textit{proper} abstract moment maps. 
\begin{Def}\label{def:cobordism}
A \textbf{cobordism} between $(M_{1},\phi_{1})$ and $(M_{2},\phi_{2})$
is an oriented orbifold $W$ and a proper abstract moment map
$\psi: W \ra \ts$, such that $\partial W = M_{1} \coprod \bar{M_{2}}$ 
and $\psi|_{M_{i}} = \phi_{i}$.
\end{Def}
The overbar indicates reversed orientation. The key to this definition
is the word ``proper''; the orbifold $W$ need not be compact,
but requiring the moment map to be proper ensures that the concept
is non-trivial. 

In the case where $(M,\phi)$ arises from an honest Hamiltonian $T$-space,
we require the cobordism to carry a presymplectic form (i.e. a 
closed 2-form) which restricts to the given forms on the boundary
components. In this case we will say that the spaces are
\textbf{cobordant as Hamiltonian $T$-spaces}.

This notion of cobordism respects reduction:
\begin{prop}[\cite{YK:1997MM}]
 \label{cob:commutewithred}
  Let $(M_{1},\phi_{1})$ and $(M_{2},\phi_{2})$ 
  be cobordant as oriented $T$-spaces with abstract moment maps. 
  Let $a \in \ts$ be a regular value of $\phi_{1}$ and $\phi_{2}$.
  Then the reduced spaces $(M_{1})_{a}$ and $(M_{2})_{a}$, with the
  orientations as defined above, are
  cobordant (as compact oriented 
  orbifolds\footnote{Note that even if $M_{1}$ and $M_{2}$ are manifolds, 
  we allow them to be cobordant by an orbifold.}). 
\end{prop}

The crucial result from \cite{GGK:1996C} and \cite{YK:1997MM}
is the localization
theorem. This can be seen as a geometric antecedent to the
theorems of Duistermaat and Heckman \cite{DH:1982} and of Jeffrey and Kirwan
\cite{JK:1995}.

From now on we will restrict to the case where $T = S^{1}$, since
this is all we need for the signature formula. This makes the 
statements simpler, although the essential picture is the same for
higher rank tori. 

We need some preliminaries about orientations.
Let $M$ be an oriented, compact $T$-space with abstract moment map
$\phi$. Denote the components of the fixed point set by 
$F_{1},\dots, F_{m}$. Assume that each component is orientable.
(In the Hamiltonian case, the fixed point components are 
in fact symplectic, so this is a reasonable assumption.) 
Fix once and for all an identification of $\tlie$ with $\R$
(a ``polarization'' of $\tlie$). Given a fixed point component
$F_{r}$, its normal bundle $N_{r}$ is an even-dimensional, oriented real
vector bundle, so it makes sense to speak of the real weights
of the action of $T$ on $N_{r}$. These are nonzero integers defined
only up to sign, and we will denote them 
$\pm \alpha_{r,k}, k=1,\dots, q_{r}$, where 
$q_{r} = \mathrm{rank} (N_{r})/2$.
By $\alpha^{\#}_{r,k}$ we will mean the weight with the
positive choice of sign; we will refer to these as polarized
weights.

Fix temporarily an orientation of $F_{r}$. This determines an
orientation of $N_{r}$, and hence makes the Pfaffian, $\Pf$,
well-defined on the fibers of $N_{r}$. Let $\xi \in \tlie$
be positive, and let $\xi_{p}$ be the isotropy action of
$\xi$ on the fiber of $N_{r}$ over some $p \in F_{r}$. Let
\[
  \varepsilon_{r} = \mathrm{sgn} \: \Pf (\xi_{p}),
\]
which is clearly independent of the choice of $p$.

If we reverse the orientation of $F_{r}$, it is clear that
the sign of $\varepsilon_{r}$ reverses as well. Hence 
there is a unique orientation of $F_{r}$ which makes 
$\varepsilon_{r}$ positive. Denote by $F_{r}^{\#}$ the
fixed point component with this choice of orientation,
which we will refer to as the polarized orientation.
(If $F_{r}$ is a single point this is just a formal sign
attached to $F_{r}$.)
We will also refer to the corresponding orientation on the
bundle as the polarized orientation.

We note that in the Hamiltonian case, this orientation does 
\textit{not} necessarily agree with the symplectic orientation
on $F_{r}$. To be precise, let $M$ be a Hamiltonian $T$-space.
Use a compatible complex structure on the normal bundle $N_{r}$ to 
identify complex weights $\alpha_{r,k}$. These are honest integers, not just
up to sign. Let $\sigma_{r}$
be the number of negative complex weights. If we use the 
symplectic orientation on $F_{r}$, it is not hard to see
that the sign of the Pfaffian of the action is 
$\varepsilon_{r} = (-1)^{\sigma_{r}}$. 
Hence the symplectic
orientation differs from the polarized orientation by 
$(-1)^{\sigma_{r}}$. In particular, when the component
$F_{r}$ is a single point $p$, the formal sign we
attach to it is exactly $(-1)^{\sigma_{r}}$.
This will be crucial in getting the correct signs
in the signature formula (\ref{cob:sigsumcircle2}).

The localization theorem says roughly that the
whole space $M$ is cobordant to the disjoint
union of the normal bundles of the fixed point compoents, with 
appropriately defined moment maps.

Fix an $N_{r}$, and equip it with an equivariant inner product.
Split the bundle up into rank two subbundles corresponding to the
weights $\alpha_{r,k}$, and denote the corresponding splitting
of any $v \in N_{r}$ by $v = \sum v_{k}$. 
Define a map on the bundle by
\[
  \psi_{r} (v) = \phi(F_{r}) + 
        \sum_{k=1}^{q_{r}} ||v_{k}||^{2} \alpha^{\#}_{r,k}.
\]
This is a proper moment map for the $T$-space $N_{r}$.

We are now ready to state the localization theorem.
\begin{thm}[\cite{GGK:1996C},\cite{YK:1997MM}]
 \label{cob:localization}
  Let $T=S^{1}$. 
  Let $M$ 
  be a compact $T$-space with abstract moment map $\phi$. 
  Fix a polarization for $T$. Define the polarized weights 
  $\alpha^{\#}_{r,k}$
  and the maps $\psi_{r}$ as above.
  Then $M$ is cobordant, as a Hamiltonian $T$-space, 
  to the disjoint union
  of the normal bundles:
\begin{equation}\label{cob:loceq}
  (M,\phi) \sim  \coprod_{r} (N_{r},\psi_{r}).     
\end{equation}
\end{thm}

(Note: with the appropriate definition of a polarization and of the
moment maps on the normal bundles, 
essentially the same statement applies in the higher-rank torus case.)

\subsection{Cobordism and Reduction}
\label{subsec:cobred}

It remains to state what happens under reduction. The localization
theorem says that the original space is cobordant to a sum
of local models, so any reduction will be cobordant to a sum of
reduced local models. In the circle case it is particularly easy to
see what these reduced local models are. 
We begin by analyzing the reduction 
of a single vector space.

Let $V$ be a real orthogonal $T$-representation with real weights 
\[
\pm \alpha_{1},\dots, \pm \alpha_{q} \ne 0.
\]
Give $V$ the polarized
orientation (so that the
Pfaffian of the action of a positive generator of $\tlie$ is positive,
as discussed above).
Fix an element
$b \in \ts$. Equip $V$ with the map $\psi_{\alpha,b}: V \ra \ts$ given by
\[
  \psi_{\alpha,b} (v) = b + 
        \sum_{k=1}^{q} ||v_{k}||^{2} \alpha^{\#}_{k}.
\]
Denote this oriented $T$-space with abstract moment map by 
$V(\alpha,b)$. As a vector space it is just $\C^{n}$.

Given such a local model $V(\alpha,b)$, and a point $a \in \ts$, $a \ne b$, 
there are two possibilities for the reduction of $V$ at $a$. 
If $a < b$ then the reduction is
empty, since the image of $\psi_{\alpha,b}$ lies to the right of $b$.
If $a > b$ then the reduction is
\begin{align*}
  V(\alpha,b)_{a} 
    &=  (\psi_{\alpha,b})^{-1}(a) / S^{1} \\
    &= \{ v \in V \: | \:  
      b + \frac{1}{2} \sum_{k=1}^{q} |v_{k}|^{2} \alpha^{\#}_{k} = a \} 
        / S^{1}  
\end{align*}
which is a twisted complex projective space of dimension $2q-2$.
Denote the (possibly empty) reduced space by $X(\alpha,a,b)$.
This is generally an orbifold, but (when it is nonempty) 
its rational cohomology ring is isomorphic to that of the ordinary
projective space $\CP^{n-1}$. 
Further, it is not hard to see that the orientation that
$X(\alpha,a,b)$ inherits as a reduced space is the same as the
complex orientation it has under the identification with
a twisted $\CP^{n-1}$, because of the condition we put
on the orientation of $V$.

The case of a vector bundle is not much more complicated.
Let $V$ be a real, orientable 
orthogonal $T$-bundle over a trivial $T$-space $F$,
with real weights 
$\pm \alpha_{1},\dots, \pm \alpha_{q} \ne 0$. Again give $E$ the
polarized orientation.
Fix an element
$b \in \ts$. Equip $E$ with the map $\psi_{\alpha,b}: E \ra \ts$ given by
\[
  \psi_{\alpha,b} (v) = b + 
        \sum_{k=1}^{q} ||v_{k}||^{2} \alpha^{\#}_{k}.
\]
This is an oriented $T$-space with abstract moment map, whose fibers
are copies of $V(\alpha,b)$.

Since $E$ fibers equivariantly over the trivial $T$-space $F$, so will the 
reduction. Since the fibers are simply copies of $V(\alpha,b)$,
the fibers of the reduced space $E_{a}$ will be copies of $X(\alpha,a,b)$.
Hence the only spaces we have to deal with are bundles whose fibers
are twisted complex projective spaces. 

Note also that $E_{a}$ is an associated bundle of the 
oriented orthonormal frame bundle
of $E$. For, the $T$-action on the fibers of $E$ commutes with the
action of the structure group,
so the structure group acts naturally on the reduction by $T$.
In particular, the structure group of the reduced space as a bundle over $F$ 
is compact and connected. 

Combining these local results with the localization theorem, and the fact that
cobordism respects reduction, we have 
\begin{thm}[\cite{YK:1997MM}]
 \label{cob:localredgeneral}
  Let $T=S^{1}$, and choose a fixed polarization of $T$.
  Let $M$ be a compact oriented $T$-space with 
  abstract moment map $\phi$. Assume that the fixed point components are
  orientable, and give them the polarized orientation $F_{r}^{\#}$. 
  Let $a \in \ts$ be a regular value of $\phi$. 
  Then there is an oriented orbifold cobordism
\[
  M_{a} \sim \coprod_{r} (N^{\#}_{r})_{a}
\]
  where each $(N^{\#}_{r})_{a}$ is a fiber bundle over $F_{r}^{\#}$
  with fiber $X_{r} = X(\alpha_{r},a,\phi(F_{r}))$, and with
  compact connected structure group. The fibers $X_{r}$ are
  orientation-preserving diffeomorphic to twisted complex 
  projective spaces. 

  If $M$ is Hamiltonian, then the polarized orientation $F_{r}^{\#}$
  differs from the symplectic orientation of $F_{r}$ by 
  $(-1)^{\sigma_{r}}$, where $\sigma_{r}$ is the number of
  negative weights of the isotropy action at $F_{r}$.
\end{thm}

\subsection{The Signature Formula in the Circle Case}
\label{subsec:sigformula}

Among the topological invariants that are preserved by cobordism, 
there is one that is singled out by being multiplicative
for any fiber bundle with oriented fiber and compact connected structure
group: the signature. We recall the 
\begin{Def}\label{cob:sigdef}
  Given a connected, compact, oriented $4k$-dimensional rational homology 
  manifold $M$, the \textbf{signature}
  $\Sign(M)$ is defined to be the signature of the symmetric bilinear
  pairing $H^{2k}(M,\Q) \tensor H^{2k}(M,\Q) \ra H^{4k}(M,\Q) \iso \Q$
  given by the cup product. If $\dim M$ is not divisible by $4$, we define
  $\Sign(M) = 0$. 
\end{Def}
Note that any orbifold is a rational homology manifold \cite{Ful:1993TV}, 
so it has a well-defined signature.

The signature is of course an important and natural invariant of
a manifold (or orbifold). 
Its special significance to this work is expressed 
by the following three results.
\begin{thm}\label{cob:sigiscobinv}
  The signature is an invariant of oriented orbifold cobordism.   
\end{thm}
\begin{proof}
  This was proved for manifolds by Thom (see Hirzebruch \cite{Hir:1966}). 
  The proof relies only on Poincar\'e duality, hence it works for
  rational homology manifolds, hence for orbifolds.
\end{proof}

The second result concerns the multiplicativity of the signature, and is a
slight modification of a result of
Chern, Hirzebruch, and Serre \cite{CheHirSer:1957}:
\begin{thm}\label{cob:chernhirzebruchserrethm}
  Let $E \ra B$ be a fiber bundle with fiber $F$ such that
  \begin{enumerate}
    \item \label{cob:basefibertotalcond} 
          $E,B,F$ are compact connected oriented orbifolds;
    \item \label{cob:structuregroupcond}
          the structure group of $E$ is compact and connected.
  \end{enumerate}
  If $E,B,F$ are oriented coherently, then 
  \[
    \Sign(E) = \Sign(B) \Sign(F).
  \]
\end{thm}
\begin{proof}
  Chern, Hirzebruch and Serre prove that the conclusion follows from
  condition (\ref{cob:basefibertotalcond}), with ``manifold'' in place
  of orbifold, and condition
  \begin{enumerate}
    \item [$2'$.] 
          The fundamental group $\pi_{1}(B)$ acts trivially on 
          $H^{*}(F)$. 
  \end{enumerate}
  As above, their proof uses only Poincar\'e duality (and the spectral
  sequence of a fibration), and hence applies equally well 
  to rational homology manifolds, such as orbifolds. Also, condition 
  (\ref{cob:structuregroupcond}) implies ($2'$).
  For, the action of $\pi_{1}(B)$ factors through the structure group; and 
  the action of a compact connected group on cohomology is trivial. 
  Our version of the theorem then follows. 
\end{proof}

To emphasize the special role of the signature, we briefly mention 
\begin{thm}[Borel, Hirzebruch; see \cite{HirBerJung:1992}]
 \label{cob:signatureisonly}
  The signature is the only rational invariant of oriented cobordism 
  that is multiplicative in the sense of Theorem 
  \ref{cob:chernhirzebruchserrethm} and takes the value 1 on
  $\CP^{2}$.
\end{thm}

Thms. \ref{cob:sigiscobinv} and \ref{cob:chernhirzebruchserrethm}, 
together with 
Thm.\ \ref{cob:localredgeneral},
immediately give the formula for the signature of the reduction by
a circle action:
\begin{align}\label{cob:localsig1}
  \Sign(M_{a}) &= \sum_{r} \Sign(F_{r}^{\#}) \Sign(X_{r}) \\
               &= \sum_{r} (-1)^{\sigma_{r}} \Sign(F_{r}) \Sign(X_{r})
\end{align}
where the second equation applies to the Hamiltonian case,
in which $F_{r}$ is given its symplectic orientation.

In the circle case the computation of 
$\Sign(X_{r})$ is quite simple, since $X_{r}$ is either empty, or a 
twisted complex projective space of complex dimension  $q_{r}-1$,
where $q_{r}$ is the half-rank of the normal bundle $N_{r}$.
The signature of a complex projective space is just $1$ when 
its complex dimension is even, $0$ when it is odd.

To be precise, fix a polarization of $T$.
Let $F_{1},\dots, F_{m}$ be the 
fixed point components, and let $\phi_{r} = \phi(F_{r})$. In the
Hamiltonian case, let
$f_{r}$ and $b_{r}$ be the number of positive and negative weights,
respectively, of $F_{r}$ (so $f_{r} + b_{r} = q_{r}$). 
Then the above formula becomes the following:
\begin{theorem}\label{thm:circlesignatureformula}
The signature of a regular reduction of a space with circle action and
abstract moment map is given by
\begin{align}
     \Sign(M_a) &= \sum_{r: \; \phi_{r} < a, \: q_{r} odd} 
                       \Sign(F_r^{\#}) \label{cob:sigsumcircle1} \\ 
                &= \sum_{r: \; \phi_{r} < a, \: q_{r} odd} 
                       (-1)^{b_r} \Sign(F_r), \label{cob:sigsumcircle2}  
\end{align}
where the latter equation applies only to the Hamiltonian case,
and $F_{r}$ is given its symplectic orientation.
\end{theorem}
This can be restated as a wall-crossing formula, which will be more
convenient later. 

\begin{figure}[htbp]
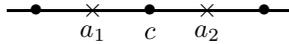

  \begin{center}
    \leavevmode
    \[ \xy/r1.8pc/:
  (0,0)="1", 
  (0.5,0)="2",
  (2.5,0)="3", "3"+/d2ex/*{c},
  (4.5,0)="4", 
  (5,0)="5", 
  (1.5,0)="6", "6"+/d2ex/*{a_{1}},
  (3.5,0)="7", "7"+/d2ex/*{a_{2}},
  "1";                     
  "5" **\dir{-},           
  "2"*{\bullet};           
  "3"*{\bullet};           
  "4"*{\bullet};           
  "6"*{\times};            
  "7"*{\times};            
\endxy \]
    \caption{Wall Crossing in the Circle Case.}
    \label{fig:circlewall}
  \end{center}
\end{figure}

Fix a singular value $c \in \ts$ of the moment map. Let
$F_{1},\dots , F_{m}$ be the fixed point components with 
$\phi(F_{r}) = c$. Pick $a_{1} < c < a_{2}$ such that there are
no other singular values in the interval $[a_{1},a_{2}]$.
(See Fig.\ \ref{fig:circlewall}.)

Denote the half-rank of the normal bundle $N_{r}$ by $q_{r}$.
Then
\begin{equation}\label{cob:wallcrossforcircleabstract}
  \Sign(M_{a_{2}}) - \Sign(M_{a_{1}}) 
     = \sum_{r=1 \dots m, \: q_{r} \mathrm{~odd}}^{m} \Sign(F_{r}^{\#}).
\end{equation}
In the Hamiltonian case, we have a slightly more explicit formula:

Define the function $w_{S}: \N \times \N \ra \Z$ by
\begin{equation}\label{cob:wsdefeq}
  w_{S}(f,b) := \left\{  
                 \begin{array}{cl}
                  (-1)^{b} & \mathrm{~if~} f+b \mathrm{~is~odd}\\
                  0        & \mathrm{~if~} f+b \mathrm{~is~even}.
                 \end{array}
                \right.
\end{equation}
Then
\begin{equation}\label{cob:wallcrossforcircle}
  \Sign(M_{a_{2}}) - \Sign(M_{a_{1}}) 
     = \sum_{r=1}^{m} w_{S}(f_{r},b_{r}) \Sign(F_{r}).
\end{equation}
Here $F_{r}$ has its symplectic orientation.

The idea of ``wall-crossing,'' and generalizations 
of (\ref{cob:wallcrossforcircle}) will be central to what follows.

\textit{Example.} \label{cob:examplecp3circle}
We first describe a whole family of simple examples of 
Hamiltonian torus actions.
Let $M = \CP^{n}$. This has a natural action of $T^{n+1}$ inherited
from the standard linear action of $T^{n+1}$ on $\C^{n+1}$ 
\[
  (\lambda_{1},\dots ,\lambda_{n+1}) \cdot (z_{1},\dots ,z_{n+1})
    = (\lambda_{1} z_{1}, \dots,  \lambda_{n+1} z_{n+1}).
\]
(This action is not effective, since the diagonal circle acts
trivially on $\CP^{n}$.) This action is easily seen to be Hamiltonian,
with moment map 
\[
  \phi([z_{1},\dots ,z_{n+1}]) 
      = \frac{1}{||z||} \sum_{k=1}^{n+1} |z_{k}|^{2} e_{k},
\]
whose image is exactly the standard $n$-simplex in $\R^{n+1}$,
\[
  \Delta_{n} = \{ a_{1},\dots, a_{n+1} \: | \: \sum a_{k} = 1 \}.
\]
Choose a subtorus
$T \subset T^{n+1}$ of rank $d$,
which does not contain the diagonal circle. Then $M$
is an effective Hamiltonian $T$-space, and the moment map is just $\phi$
composed with the projection $\R^{n+1} \ra \ts$. The moment polytope
is the image of the standard $n$-simplex under this projection.
If $T$ is generic, its fixed points are the same as those
of $T^{n+1}$, which are clearly the $n+1$ isolated points
$[z_{1},\dots ,z_{n+1}]$ with $z_{j} = \delta^{j}_{k}$, for
$k=1,\dots ,n+1$.

\begin{figure}[htbp]
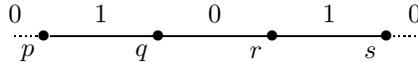

  \begin{center}
    \leavevmode
    \[ \xy/r1.8pc/:
  (0,0)="1", "1"+/u2ex/*{0}, 
  (0.5,0)="2", "2"+/dl2ex/*{p}, 
  (2.5,0)="3", "3"+/dl2ex/*{q}, 
  (4.5,0)="4", "4"+/dl2ex/*{r}, 
  (6.5,0)="5", "5"+/dl2ex/*{s}, 
  (7,0)="6", "6"+/u2ex/*{0}, 
  (1.5,0)="7", "7"+/u2ex/*{1}, 
  (3.5,0)="8", "8"+/u2ex/*{0}, 
  (5.5,0)="9", "9"+/u2ex/*{1}, 
  "1";                     
  "6" **\dir{.},           
  "2"*{\bullet};           
  "5" **\dir{-},           
  "3"*{\bullet};           
  "4"*{\bullet};           
  "5"*{\bullet};           
\endxy \]
    \caption{Wall Crossing Example: $\CP^{3}$ with Circle Action.}
    \label{fig:cp3circle}
  \end{center}
\end{figure}

For now, choose $n=3$, $d=1$, so we get a generic Hamiltonian
circle action on $\CP^{3}$, with $4$ isolated fixed points.
The moment polytope is shown in Figure \ref{fig:cp3circle},
with the images of the fixed points shown as $p,q,r,s$.
Since the complex dimension is $3$, each fixed point
has $3$ weights associated to it. It is not hard to see
by looking at the action that at $p$, all three weights
point to the right; at $s$, all three point left; at
$q$, one points left, two point right; and at $r$,
two point left and one points right.

Since all fixed points are isolated,
the wall-crossing rule (\ref{cob:wallcrossforcircle}) gives
the signatures shown in the figure. These agree with a direct
calculation of the reduced spaces, which easily shows
that the reductions in the outer two chambers are 
(twisted) $\CP^{2}$, which has signature $1$, 
and the reduction in the inner
chamber has the rational cohomology of $\CP^{1} \times \CP^{1}$, 
which has signature $0$.

\subsection{Singular Reduction}\label{subsec:singularreduction}

Recently Lerman and Tolman \cite{LermanTolman:1998Sing} 
have found a beautiful way
to analyze the rational intersection cohomology ring of singular reductions
for Hamiltonian circle actions. We can use a piece of their
result to easily extend the above formula to the
singular case. The key is that Thm.\ \ref{thm:circlesignatureformula} 
applies to any abstract moment map. 

Given a singular value $c$ of the moment map, let 
$\F_{c} = \{ F_{1},\dots, F_{m} \}$
be the set of fixed point components with $\phi(F_{k}) = c$, as above.
Divide $\F_{c}$ into two disjoint sets: 
\begin{align*}
  \F_{c}^{+} &= \{ F_{k} \: | \: b_{k} < f_{k} \} \\
  \F_{c}^{-} &= \{ F_{k} \: | \: b_{k} \ge f_{k} \}.
\end{align*}

Lerman and Tolman show that by deforming the moment map $\phi$ 
to a new map $\tilde{\phi}$, which is now only an abstract moment map,
they can find a ``small resolution'' of the singularities of 
the singular reduction $\phi^{-1}(c)/S^{1}$. Precisely,

\begin{theorem}[Lerman--Tolman]
\label{thm:LermanTolmanSmallResolution}
  Let $(M,\omega,\phi)$ be a Hamiltonian $S^{1}$-space.
  Let $c$ be a singular value of $\phi$, and define 
  $F_{1},\dots ,F_{m}$, $\F_{c}^{\pm}$ as above. 
  Then there exist $a,b \in \R$, $a < c < b$, and 
  an abstract moment map $\tilde{\phi}$, such that: 
\begin{enumerate}
\item $a,b$ are regular values of both $\phi$ and $\tilde{\phi}$;
\item there are no fixed points in $\phi^{-1}([a,b])$ besides those
      in $\F_{c}$;
\item $\tilde{\phi}$ and $\phi$ are equal on $M \setminus \phi^{-1}([a,b])$;
\item $c$ is a regular value of $\tilde{\phi}$;
\item $a < \tilde{\phi}(F_{k}) < c$ for all $F_{k} \in \F_{c}^{-}$
  and $c < \tilde{\phi}(F_{k}) < b$ for all $F_{k} \in \F_{c}^{+}$;
\item $\tilde{M}_{c} := \tilde{\phi}^{-1}(c)/S^{1}$ 
  is a small resolution
  of $M_{c} = \phi^{-1}(c)/S^{1}$. In particular there is a
  pairing-preserving isomorphism between the
  (middle perversity) intersection cohomology of $M_{c}$
  and the ordinary cohomology of $\tilde{M}_{c}$.
\end{enumerate}
\end{theorem}

We can now apply our formulas to the new abstract moment map
$\tilde{\phi}$ to get a formula for the singular reduction.
We simply compare $\tilde{M}_{c}$ to $M_{a}$ (or to $M_{b}$).

The wall-crossing formula (\ref{cob:wallcrossforcircle}) 
applied to $\tilde{\phi}$ gives 
\[
  \Sign(\tilde{M}_{c}) 
    = \Sign(\tilde{M}_{a}) + 
        \sum_{F_{r} \in \F_{c}^{-}, \quad  q_{r} \mathrm{~odd}} \Sign(F_{r}^{\#}).
\]
One gets a similar formula involving $\F_{c}^{+}$ by comparing
to $\tilde{M}_{b}$. 

Note that since the $F_{r}$ are actually fixed point components of the
original Hamiltonian action, they are symplectic, and the
normal bundles do have complex structures. So we can
restate this formula in terms of the symplectic orientation of $F_{r}$
and the number of negative weights, just as in the honest Hamiltonian
case:
\[
  \Sign(\tilde{M}_{c}) 
    = \Sign(\tilde{M}_{a}) + 
        \sum_{F_{r} \in \F_{c}^{-}} w_{S}(f_{r},b_{r}) \Sign(F_{r}).
\]

Since $\tilde{M}_{c}$ is a small resolution of $M_{c}$, their
signatures are equal (where the signature of the singular space
$M_{c}$ is understood to be taken in intersection cohomology).
Also, since $\tilde{\phi}$ and $\phi$ only differ near the 
level $c$, we have $\tilde{M}_{a} = M_{a}$, and $\tilde{M}_{b} = M_{b}$.

Hence we have proved
\begin{theorem}\label{thm:singularsignature}
  Let $(M,\omega,\phi)$ be a compact Hamiltonian $S^{1}$-space,
  and let $c$ be a singular value of the moment map. Define
  the sets $\F_{c}^{\pm}$ and the function $w_{S}$ as above.
  Choose $a,b \in \R$, $a < c < b$, such that there are no critical
  values of $\phi$ in the intervals $[a,c)$ and $(c,b]$. 
  Then 
\begin{align}\label{singularsignatureformula}
  \Sign(M_{c}) 
    &= \Sign(M_{a}) + 
        \sum_{F_{r} \in \F_{c}^{-}} w_{S} (b_{r},f_{r}) \Sign(F_{r}) \\
    &= \Sign(M_{b}) - 
        \sum_{F_{r} \in \F_{c}^{+}} w_{S} (b_{r},f_{r}) \Sign(F_{r}) 
\end{align}
\end{theorem}

Note that the jump between the signatures of the regular reductions 
$\Sign(M_{a})$ and $\Sign(M_{b})$, which is equal to 
\[
  \sum_{F_{r} \in \F_{c}} w_{S} (b_{r},f_{r}) \Sign(F_{r}) 
\]
is simply divided into two steps, one between 
$\Sign(M_{a})$ and $\Sign(M_{c})$, the other between 
$\Sign(M_{c})$ and $\Sign(M_{b})$. Each
fixed point component contributes to one step or the other 
depending on whether it has more weights up or down.

One may wonder why in the definitions of $\F_{c}^{\pm}$ we chose
the $F_{r}$ with $b_{r} = f_{r}$ to lie in $\F_{c}^{-}$
and not in $\F_{c}^{+}$. This is in fact an arbitrary choice;
Lerman and Tolman's theorem is still true if we put these
$F_{r}$ in $\F_{c}^{+}$ instead. Precisely, making this choice
gives a different small resolution ${\tilde{M}'}_{c}$ of $M_{c}$,
with a cohomology ring that is not isomorphic to the one resulting
from our choice; however, the pairings are isomorphic. 
This makes sense in Thm.\ \ref{thm:singularsignature}, 
because the fixed point components with $b_{r} = f_{r}$
do not contribute to the wall-crossing:
when $b_{r} = f_{r}$, $b_{r} + f_{r}$ is even and $w_{S}(b_{r},f_{r})=0$.
Hence it doesn't matter where we include them.


\section{The Structure of the Moment Polytope}
\label{xraychap}

To describe the formula for the signature in the torus case,
we first need to review what is known about the structure of the moment
polytope, and introduce some notation.

Let $T$ be a torus of rank $d$. Let $(M,\omega,\phi)$ be a compact 
Hamiltonian $T$-space of dimension $2n$. 
For convenience, assume the action is effective. 
By the convexity theorem of 
Atiyah and Guillemin-Sternberg \cite{A:1982CCH} \cite{GS:1982C}, 
the moment image $\Delta = \phi(M)$ is a convex 
polytope, and is in fact the convex hull of the image the fixed point
set $M^{T}$. However, one can say a great deal more about the structure
of $\Delta$. We will introduce a refinement of Tolman's notion
of the \textit{X-ray} of $M$ to provide the framework for later 
calculations. 

Note: for the remainder of this section, we will refer to
connected subgroups of $T$ as subtori, to distinguish them from
arbitrary subgroups. Since the results we are interested in 
are rational phenomena (as opposed to torsion information) and hence are not
affected by discrete stabilizers, we will ignore disconnected subgroups.

Recall from the structure theory of compact transformation groups
(see e.g. \cite{KK:1991}) that $M$ decomposes into a finite
set of orbit-type strata. Given a subgroup $H \subset T$, 
the orbit-type stratum $M_{H}$ is the set of points $p \in M$
such that the isoptropy subgroup $T_{p} = H$. Let $X_{1},\dots, X_{m}$
be the set of connected components of orbit-type strata, with 
corresponding isotropy subroups $T_{1},\dots ,T_{m}$. For each
$X_{j}$ denote the closure by $F_{j} = \bar{X_{j}}$. We will
abuse notation and refer to these as the \textit{strata} of $M$. 
Let $\F = \{ F_{1},\dots ,F_{m} \}$ be the set of strata. Each $F_{j}$
is a connected component of the fixed point set $M^{T_{j}}$
of the subtorus $T_{j}$\footnote{Not 
all connected components of $M^{T_{j}}$ are of this form;
for example a connected component of $M^{T}$ will often be
a connected component of $M^{T_{j}}$ as well, but it is not
the closure of a stratum corresponding to $T_{j}$.}.
Note that $T_{j}$ is the stabilizer of a generic point on $F_{j}$, but
some points in $F_{j}$ will have larger stabilizer.
Denote  $\dim(T_{j})$ by $d_{j}$. 

Each $F_{j}$ is a symplectic 
manifold in its own right, by the equivariant Darboux theorem 
\cite{GS:1984STP}. In fact, the restriction of the moment map
$\phi$ makes $F_{j}$ into a Hamiltonian $T$-space. However, the $T$-action
on $F_{j}$ is clearly not effective. $F_{j}$ inherits an effective action of 
the quotient torus 
$H_{j} = T/T_{j}$. The $H_{j}$-action is Hamiltonian, but the
moment map is unique only up to addition of a constant.
It is specified as follows. Consider the three exact sequences
\begin{equation}
\label{x:exactforfixed}
  \xymatrix{ 1 \ar[r] & {T_{j}} \ar[r] & {T} \ar[r] &
             {H_{j}} \ar[r] & 1 \\
             0 \ar[r] & {\tlie_{j}} \ar[r] & {\tlie} \ar[r] &
             {\h_{j}} \ar[r] & 0 \\
             0 & {\ts_{j}} \ar[l] & {\ts} \ar[l] &
             {\hs_{j}} \ar[l] & 0 \ar[l] \\
               & d_{j} & d & (d - d_{j}) &  
           }
\end{equation}
(the last row notes the respective dimensions). 
Note $\hs_{j} = \Ann(\tlie_{j})$. The moment map 
$\phi : M \ra \ts$ restricted to $F_{j}$ lies in a translate
of $\hs_{j} \subset \ts$ (since the vector fields generated by
$\tlie_{j}$ vanish on $F_{j}$). In particular, if $T_{j} = T$
then $\phi(F_{j})$ is a point. The restricted moment
map, shifted by a constant to land in $\hs_{j}$, is a moment map
for the action of $H_{j}$. 

Hence, by the convexity theorem, the image 
$\phi (F_{j})$ 
is a convex polytope in its own right, though of dimension $d-d_{j}$.
We will call these polytopes \textit{walls}.
It is important to note that two different $F_{j}$ can
have overlapping, or even identical, images under $\phi$. 
Tolman \cite{ST:1996NK} introduced the notion of the \textit{X-ray}
of $M$ to keep track of this data. 
\begin{Def}\label{def:xray}
  The \textbf{X-ray} of a compact Hamiltonian $T$-space as above is the
  family $\{ \phi(F_{j}) \}$ of convex polytopes in $\ts $,
  indexed by the set $\F$, where $\F$ is considered as a partial 
  order under inclusion.
\end{Def}
More technically, the data of the X-ray are the partial order $\F$ and the map
$F_{j} \mapsto \phi(F_{j})$, which we will also denote by $\phi$. 
We will often refer to an X-ray by the pair $(\F,\phi)$.
Note that $F_{j} \subset F_{k}$ implies 
$\phi(F_{j}) \subset \phi(F_{k})$ but not vice versa. 
We will refer to the walls coming from $T$-fixed point components
(which are just single points of $\ts$) 
as the \textit{vertices} of the X-ray.

X-rays have some nice formal properties (which are explored further
in \cite{Metzler:1997}). 
Denote the affine
span of a set $S$ in $\ts$ by $\Aff(S)$, and the linear span of $S$
(the unique linear subspace parallel to $\Aff(S)$) by $\Lin(S)$.
\begin{prop}\label{x:xrayproperties}
  Given a Hamiltonian $T$-space $M$, its X-ray map $F_{j} \mapsto \phi(F_{j})$
  is an order-preserving map from $\F$ to the set of convex polytopes in
  $\ts$, satisfying: 
  \begin{enumerate}
  \item \label{x:face} Given $F_{j} \in \F$, and a face $\delta$ of 
        $\phi(F_{j})$,
	there is a unique $F_{k} \le F_{j}$ such that $\phi(F_{k}) = \delta$.
  \item \label{x:uniquegoingup} Given $F_{j} \in \F$, and $F_{k} \ge F_{j}$, 
        $F_l \ge F_{j}$, with $\Lin(\phi(F_{k})) = \Lin(\phi(F_l))$, then 
        $F_{k} = F_l$.
  \item \label{x:dimdiff} 
        $F_{k} < F_{j} \implies \dim \phi(F_{k}) < \dim \phi(F_{j})$. 
  \end{enumerate}
\end{prop}
\begin{proof}
  First, note that $\Lin(\phi(F_{k})) = \Ann(\tlie_{k})$
  by the definition of a moment map. 
  Hence condition \ref{x:dimdiff} in the definition
  of an X-ray follows from the fact that for any $F_{j} \subsetneq F_{k}$,
  $T_{j} \supsetneq T_{k}$, and hence $\tlie_{j} \supsetneq \tlie_{k}$, 
  since $T_{j}, T_{k}$ are connected. Also, condition \ref{x:uniquegoingup} 
  is clear: the condition $\Lin(\phi(F_k)) = \Lin(\phi(F_l))$ implies
  that $F_k$ and $F_l$ are connected components of the fixed point set
  of the same subtorus. But they both contain $F$, hence they must
  be equal.

  Condition \ref{x:face} is a bit deeper. It comes from the result
  of Atiyah \cite{A:1982CCH} that the inverse
  image of any point under the moment map is connected. 
  For, consider a stratum $F_{j}$. This is a Hamiltonian
  $T$-space. Given a face $\delta$ of $\phi(F_{j})$, 
  let $\h = \Ann(\Lin(\delta))$, and let $H \subset T$ be the corresponding
  Lie subgroup. (This exists since the polytope $\phi(F_{j})$ is 
  rational \cite{G:1994MM}\footnote{This is another 
  feature of Hamiltonian $T$-spaces that will
  not be particularly relevant in this paper.}.)
  Consider the moment map of the $H$-action on $F_{j}$; this is just
  the projection of the moment map by 
  $\pi_{H}: \ts \ra \hs = \ts/\Lin(\delta)$.
  This moment map has an extremum at $\pi_{H}(\delta)$,
  which is certainly a critical value. Moreover, any point in $F_{j}$ mapping
  to this point will necessarily be a critical point. Hence the entire 
  inverse image of $\pi_{H}(\delta)$ under the $H$-moment map, 
  which is the inverse image
  of $\delta$ under the $T$-moment map, is fixed by the action of $H$.
  Moreover, this set is connected by Atiyah's result.
  Hence it is exactly some $F_{k}$, and 
  $H = T_{k}$, $\phi(F_{k}) = \delta$, and no other $F_{l} \subset F_{j}$ 
  can map to $\delta$.
\end{proof}

\textit{Remark.\/} 
Our definition is not completely consistent with
Tolman's definition of an X-ray. 
In defining the X-ray of a Hamiltonian action, 
she considers all possible stabilizer
groups, not just connected ones, and hence her X-rays
do not satisfy condition \ref{x:dimdiff} above. 
By considering disconnected
stabilizers she gets at torsion information, in particular, the 
intrinsic stabilizers of singular points in the reduced spaces
(when these spaces are orbifolds). However, we are primarily interested
in rational invariants, which are insensitive to this information. 
We could refer to our X-rays as ``X-rays modulo torsion.''
Note, however, that much of what we say in the rest of this section 
applies in a slightly modified form to Tolman's X-rays.

\begin{figure}[htbp]
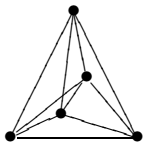

  \begin{center}
    \leavevmode
    \[ \xy/r1pc/:
  (0,0)="1",
  (4,0)="2", 
  (2,4)="3", 
  (1.6,0.75)="4", 
  (2.4,1.9)="5", 
  "1"*{\bullet};                        
  "2" **\dir{-}, "2" *{\bullet},        
  "3" **\dir{-}, "3" *{\bullet},        
  "4" **\dir{-}, "4" *{\bullet},        
  "5" **\dir{-}; "5" *{\bullet};        
  "2" **\dir{-},                        
  "3" **\dir{-},                        
  "4" **\dir{-};                        
  "2" **\dir{-},                        
  "3" **\dir{-};                        
  "2" **\dir{-};                        
\endxy \]
    \caption{X-ray of $M = \CP^4$, with 2-torus action.}
    \label{fig:cp4}
  \end{center}
\end{figure}

\textit{Example.} 
We again consider a complex projective space, as on page 
\pageref{cob:examplecp3circle}, 
but now we take $n=4$, $d=2$\footnote{We choose $d=2$ for later purposes, 
so that the signatures of the reduced spaces are nonzero. 
See Section \ref{recursivechap}.}.
For a generic choice
of the subtorus, the image of the moment polytope will appear as in
Figure \ref{fig:cp4}. (This picture is in $\R^{2}$.)
The five dots are the vertices (images of the $T$-fixed point components),
while the lines are images of strata corresponding to circles in $T$.

\subsection{Weight Data}\label{subsec:weightdata}

We also want to  record information about the weights of the
infinitesimal action of $T$ on the normal bundles to the strata.
This data will appear in the wall-crossing formulas for the signature
and other invariants. 
Moreover, the weight data at the $T$-fixed points 
often completely determine the X-ray (see \cite{Metzler:1997}).

Recall from the introduction 
that the $T$-invariant symplectic structure on $M$ 
defines a natural isotopy class of compatible, 
invariant almost complex structures.
Hence, given any point
$p \in M^{T_{j}}$, the complex weights of the action of $T_{j}$ on $T_{p} M$
are well-defined as vectors in $\ts$.
Given a stratum $F_{j} \subset M^{T_{j}}$, 
denote the weights of the $T_{j}$ action on $T_{p} M$, for any 
$p \in M$, by $\alpha_{j}= \{ \alpha_{j,k} \in \ts_{j} \}$. 
(Since $F_{j}$ is connected this is independent of $p$.)
\begin{Def}\label{x:newweightedxraydef} 
  The \textit{weighted X-ray} of $M$ is the X-ray together with
  the assignment to each stratum $F_{j}$ of the weights
  $\alpha_{j,k} \in \ts_{j}$. 
\end{Def}
Note that we consider all of the weights on 
$TM|_{F_{j}}$, not just the weights of the normal bundle $N_{F_{j}}$.
Of course the additional weights of $T_{j}$, corresponding to $T F_{j}$, are
zero, but we will see below that the bookkeeping works out better if we 
include these weights as well. Also, $\alpha_{F}$ is not a set of vectors, but
a family, as we need to keep track of multiplicities.

We need to describe some of the properties that weighted X-rays
of Hamiltonian $T$-spaces have. The most important is a
consequence of the equivariant Darboux theorem, which says that
the weights determine the X-ray ``locally.''

We will use the
following definition of a \textit{local model}. In any vector space $A$,
given a set of vectors $\alpha = \{\alpha_{1},\dots,\alpha_{n} \in A \}$, 
denote the positive cone generated by $\alpha$ by
\begin{equation}\label{coneeq}
     \Cone(\alpha_{1},\alpha_{2},\dots \alpha_{n}) 
        = \left\{ \left. \sum_{k=1}^{n} t_{k} \alpha_{k} 
                   \: \right| \: t_{k} \ge 0 \right\}. 
\end{equation}

Call a subset
$S = \{ \alpha_{i_{1}},\alpha_{i_{2}},\dots \alpha_{i_{j}}\} \subset \alpha$ 
\textit{linear} if no other $\alpha_{i} \in \alpha$, $\alpha_{i} \notin S$ 
lies in the linear span of $S$;  
in other words, if $S = \alpha \intersect V$ for some linear subspace 
$V \subset A$.
\begin{Def}\label{x:localmodeldef}   
   The \textbf{local model} generated by 
   $\{\alpha_{1},\dots,\alpha_{n} \in A \}$ is the set of cones 
     \[
        \{\Cone (S) \: | \: S \mathrm{~a~linear~subset~of~} \alpha  \}.
     \]
\end{Def}
We will also refer to the intersection of these cones with any neighborhood
$U$ of $0$ as a local model.
 
We will also need to consider the natural
maps $  \pi_{jk}: \ts_{j} \ra \ts_{k}$ arising when
$F_{j} \subset F_{k}$, and the quotient map 
$\pi_{j}: \ts \ra \ts_{j}$.

\begin{prop}\label{x:weightedxraydef}
   The weighted X-ray $(\F,\phi,\alpha)$ of a Hamiltonian $T$-space
   satisfies:
   \begin{enumerate}
     \item  
       \label{x:darbouxcond}
       for every stratum $F_{j}$,
       there is a neighborhood $U \subset \ts$ of $\phi(F_{j})$ such that
       $\{\pi_{k}(\phi(F_{k}) \intersect U) \: | \: F_{k} \ge F_{j} \}$ 
       is the local model generated by $\alpha_{j}$;
     \item 
       \label{x:consistcond}
       for every $F_{j}, F_{k} \in \F$, $F_{j} \le F_{k}$, 
       the families of vectors 
       $\pi_{jk}(\alpha_{j})$ and $\alpha_{k}$ are equal, up to 
       rearrangement. (I.e. they contain the same vectors with the
       same multiplicities.)
\end{enumerate}
\end{prop}
We will refer to condition \ref{x:darbouxcond} as the Darboux 
condition, and condition \ref{x:consistcond} as the consistency condition. 
\begin{proof}
  The weights of the action of $T_{k}$ on the tangent bundle 
  can be determined by looking at the tangent space to any point 
  $p \in F_{k}$. Suppose we have $F_{j} \subset F_{k}$. Choose the
  point $p$ to be in $F_{j}$. Then the representation of $T_{j}$ on
  the tangent space $T_{p}(M)$ must restrict to the representation
  of $T_{k}$. In terms of weights, this is exactly condition 
  \ref{x:consistcond} in the proposition.

  To show that the Darboux condition holds for an arbitrary wall,
  first note that it is enough to show that it holds for the
  fixed points of the whole torus. For, given a stratum
  $F_{j} \subset M^{T_{j}}$, we can restrict the action to $T_{j}$; the 
  effect on the moment map is to compose with the projection
  $\pi_{j}:\ts \ra \ts_{j}$. Clearly the formation of a local model, and
  hence the Darboux condition, respects
  this projection, so $F_{j}$ will satisfy the Darboux condition
  with respect to the $T$ action iff it satisfies it with respect to
  the $T_{j}$-action.

  So assume that $p \in F \subset M^{T}$. Denote the
  weights of the action of $T$ on $T_{p}M$ by 
  $\alpha = \{ \alpha_{1},\dots \alpha_{n} \}$.
  The Darboux theorem, in the equivariant setting, says
  (\cite{GS:1984STP}, p. 251) that we can equivariantly identify
  a neighborhood of $p$ with $(\C^{n},\omega_{std},\phi_{\alpha})$,
  where $\C^{n}$ has the action of $T$ given by the weights $\alpha_{k}$,
  that is, the infinitesimal action of $v \in \tlie$ is given by
  \[
    v \cdot z = (i \alpha_{1}(v) z_{1}, \dots, i \alpha_{n}(v) z_{n}), 
  \]
and
  \[
    \omega_{std} = \frac{i}{2} \sum_{k=1}^{n} dz \wedge d\bar{z},
  \]
  \[
    \phi_{\alpha}(z) = \frac{1}{2} \sum_{k=1}^{n} |z_{k}|^{2} \alpha_{k}.
  \]  
  Hence we only need to determine the relationship of the weights
  to the images of the fixed point sets of subtori for this linear
  space. 

  Given a subspace $W \subset \tlie$, let $V = \Ann(W)$.
  The set of points in $\C^{n}$ infinitesimally fixed by all $w \in W$
  is 
\[
  (\C^{n})^{W} = \{(z_{1},\dots ,z_{n}) \: | \: 
     \forall k, z_{k} \ne 0 \implies \alpha_{k} \in \Ann(W)\}
\]
  and its moment image is 
\begin{align}\label{x:darbouxcomputeeq}
  \phi((\C^{n})^{W}) &= \left\{ \left. 
         \frac{1}{2} \sum_{k=1}^{n} |z_{k}|^{2} \alpha_{k} \: \right| \:
                         z_{k} \ne 0 \implies \alpha_{k} \in V \right\} \\
                     &= \Cone(\alpha \intersect V). 
\end{align}
  Since any $V$ is the annihilator of some $W$, the walls 
  are exactly the sets of the form (\ref{x:darbouxcomputeeq}), which
  proves the Darboux condition.
  
\end{proof}

The Darboux condition involves both existence and uniqueness: 
for every linear subset $S$ of the weights of $F_{j}$, 
there is a unique wall, comparable to $F_{j}$, that locally looks like
the cone on $S$. However, the uniqueness only applies to strata
$F_{k} \ge F_{j}$; there may be other walls which happen to lie near
$F_{j}$ (or even overlap it) but if they are not comparable to $F_{j}$
they need have no particular relation to it. 

Note that one consequence of the Darboux condition is that the
weights of the $T$-fixed point components 
lie along the one-dimensional walls of the X-ray.
Since for our purposes we will not need to know the length of the weights, 
it is not always necessary to draw the weights explicitly, since we
can read off their directions from the 1-skeleton of the X-ray.
However, there are often multiple weights pointing in the same
direction, and we will need to be careful to record this multiplicity
information.

In the example of $\CP^{4}$ above (Fig. \ref{fig:cp4}),
each $T$-fixed point is isolated, and hence has four weights
attached to it. Clearly there must be exactly one weight lying in 
the direction of each line emanating from a given vertex. 
Here there is no ambiguity in the multiplicities, and we can
read off all of the information we will need from the figure.

So far we have concentrated on the walls of the X-ray, which 
come from the various fixed point sets. Now we 
look at the complement of the walls. Assume from now on that $M$ is connected.
It is easy to see that the 
union of the walls (not counting the largest ``wall,'' 
the entire moment polytope) is the set of singular values of the moment map,
and hence the complement is the set of regular values:
\begin{equation}\label{x:deltaregdef}
  \Delta_{reg} = \Delta \setminus \Union_{F \ne M} \phi(F).
\end{equation}
As mentioned in the introduction, 
this is an open set with a finite number of components, which
we call \textit{chambers}. It is well-known that these chambers
are open convex polytopes. (In fact the convexity follows from Props. 
\ref{x:xrayproperties} and \ref{x:weightedxraydef}, as we prove
in \cite{Metzler:1997}.)

We will want the following more general definition when we discuss
recursive invariants in Chap.\ \ref{recursivechap}.
Recall that a stratum $F_{j} \subset M^{T_{j}}$
is an effective Hamiltonian $T/T^{j}$ space. 
The regular values of the restricted moment map are clearly
\[
  \Reg(F_{j}) = \phi(F_{j}) \setminus \Union_{F_{k} < F_{j}} \phi(F_{k}),
\]
which is a relatively open set in $\phi(W)$, again a finite union of
(relatively) open convex polytopes \cite{Metzler:1997}.
We need to be careful to remember which $F_{j}$ these
came from, so call a pair $(F_{j},P)$, where $P$ is a component of 
$\Reg(F_{j})$, a \textit{subchamber} of the X-ray. 
By abuse of notation we will often
refer to $P$ as a subchamber, when it will cause no confusion.
Note that the vertices of the X-ray are subchambers: the
corresponding torus $T/T$ is trivial, so its Lie algebra is zero-dimensional,
and every point is a regular point for the ``moment map.''


\section{Recursive Invariants of X-rays}
\label{recursivechap}

Given two adjacent chambers in the moment polytope of
a Hamiltonian $T$-space, we would like to have a wall-crossing
formula that relates the signature of a reduction on one
side to that of a reduction on the other side. The simplest
analogy with the circle case would be if the difference were expressed
in terms of the weights associated to the wall, and 
the signature of some fixed point component(s). 
The truth is slightly more complicated. In fact the signatures
of the reduced spaces are best understood as part of a larger
system, which we will refer to as a \textit{recursive invariant}.

We know from Prop.\ \ref{invarianceofreduction} 
that topological invariants of reductions, like the signature, 
are functions only of the set of chambers of the moment polytope.
We want to elaborate on this idea slightly. So far, we have only
considered the set of regular values $\Delta_{reg}$ as places
where we can reduce. However, every point in the moment polytope 
is a regular value for some Hamiltonian action. Precisely, let
$F_{j} \in \F$ be a stratum corresponding to the subtorus $T_{j}$,
and fix a subchamber $P \subset \phi(F_{j})$.
For any point $q \in P$, $q$ is a regular value of the moment map
for the action of the quotient torus $H_{j}$\footnote{We must technically 
pick a fixed translation of $\Aff(\Delta_{j})$ into $\Ann(\tlie_{j})$ to
say this, but the formula (\ref{subreductioneq}) doesn't depend on 
this choice.}, so we can form the regular symplectic reduction
\begin{equation}\label{subreductioneq}
  (F_{j}/ \! \! /H_{j})_{q} = (\phi^{-1}(q) \intersect F_{j})/H_{j}
       = (\phi^{-1}(q) \intersect F_{j})/T.
\end{equation}
We will denote this by $M^{j}_{(q)}$, and when there is no ambiguity,
simply $M_{(q)}$. 
Note that if $q$ lies in $\Delta_{\mathrm{reg}}$ this is simply
the usual symplectic reduction of $M$ at $q$. However, when 
$q \notin \Delta_{\mathrm{reg}}$ it is to be carefully distinguished 
from the singular reduction
$\phi^{-1}(q)/T$, which is not even an orbifold in general. (In fact
it is a symplectic stratified space \cite{LerSja:1991}, 
and $M^{j}_{(q)}$ is one of the strata.)
We must make the usual caveat that this space depends not only on
$q \in \ts$ but also on the fixed point component $F_{j}$. If there
are two $F_{j},F_{k}$ whose moment images both contain $q$, we must
specify which fixed point component we are restricting to before
taking the reduction.

The recursive formula for the signature, Thm.\ \ref{rec:sigwallthm}, 
in fact calculates the 
signature on all of these ``subreductions.'' Clearly,
however, the bookkeeping threatens to get a little involved. 
The weighted X-ray summarizes exactly the data about the Hamiltonian
action that is needed.

Fix a Hamiltonian $T$-space $(M,\omega,\phi)$. Denote the
weighted X-ray of $M$ by $(\F,\phi,\alpha)$.
Let the set of all subchambers $(F_{j},P)$ of the X-ray be denoted by
$Q_{\F}$. Define a function $S: Q_{\F} \ra \Z$ by
\[
 S(F_{j},P) = \Sign(M^{j}_{q})
\]
for any $q \in P$. We will show that this function satisfies
a simple recursive formula. To make things general enough to
accomodate the Euler characteristic and the Poincar\'e polynomial
as well, we introduce a general definition. First we need a lemma,
concerning how subchambers of an X-ray can meet. 
Given a wall $\phi(F)$ of an X-ray, call a codimension 1
subwall of $\phi(F)$ a \textit{principal} subwall.

\begin{lemma}\label{rec:chamberoverlaplemma}
  Let $(\F,\phi,\alpha)$ be the weighted X-ray of $(M,\omega,\phi)$.
  Let $F \in \F$, and let $(F,P_{1})$, $(F,P_{2})$ be subchambers 
  whose closures intersect in a set which has codimension 1 in $\phi(F)$. 
  Clearly $P_{12} := \bar{P_{1}} \intersect \bar{P_{2}}$ is a convex polytope.
  Call its affine span $S$.
  Then for every principal subwall $G < F$ such that $\phi(G) \subset S$ and 
  $\phi(G) \intersect P_{12} \ne \emptyset$, 
  $P_{12}$ lies entirely in one subchamber $(G,R)$ of $G$.
\end{lemma}

\begin{figure}[htbp]
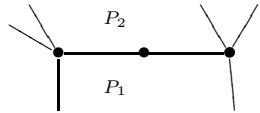

  \begin{center}
    \leavevmode
    \[ \xy/r1.8pc/:
  (0,0)="1", 
  (3,0)="2", 
  (0,-1)="3", 
  (3.1,-1)="4", 
  (-.85,0.5)="5", 
  (-0.5,0.85)="6", 
  (2.5,0.85)="7", 
  (3.5,0.85)="8",
  (1,-0.6)="9", "9"*{\scriptstyle P_1},
  (1,0.6)="10", "10"*{\scriptstyle P_2},
  (1.5,0)="11",
  (4.2,-0.6)="15", 
  (4.2,0.6)="16",  
  "1"*{\bullet};           
  "3" **\dir{-},           
  "5" **\dir{-},           
  "6" **\dir{-},           
  "2"*{\bullet} **\dir{-}; 
  "2" *{\bullet};
  "4" **\dir{-},           
  "7" **\dir{-},           
  "8" **\dir{-},           
  "11" *{\bullet}          
\endxy \]
    \caption{Disallowed Meeting of Chambers.}
    \label{fig:badchamber}
  \end{center}
\end{figure}

See Fig.\ \ref{fig:badchamber} for a picture of a situation
that this lemma rules out, where there are two subchambers
in the same wall separating $P_{1}$ and $P_{2}$.

\begin{proof}
  Since $\phi(F)$ is an X-ray in its own right, 
  it is no loss of generality to assume that $M$ is connected, that
  $F = M$, and $P_{1}, P_{2}$ are
  chambers, intersecting in a principal wall. 

  Assume that $\bar{P_{1}}$ meets $\bar{P_{2}}$ in two or more 
  (principal) subchambers of $G$. Two of them must be separated by 
  a principal subwall $W$ of $G$. Now $G$ must have weights pointing in
  the directions of both $P_{1}$ and $P_{2}$. 
  Hence by the consistency condition,
  $W$ must have weights $\beta_{1}, \beta_{2}$ 
  pointing out of $G$ in the directions
  of $P_{1}$ and $P_{2}$, respectively. 
  But by the Darboux condition there is a principal subwall
  of $F$ generated by the codimension $2$ wall $W$ and the weight 
  $\beta_{1}$, and likewise for $\beta_{2}$; these walls would break
  up $P_{1}$ and $P_{2}$, giving a contradiction.  
\end{proof}

We will be looking at the following situation.  Consider a stratum
$F \in \F$, and two subchambers $(F,P_{1}),(F,P_{2})$
separated by principal subwalls $\phi(G_{1}),\dots ,\phi(G_{m})$, 
which have a common affine span $S$.
By the lemma, in each wall $\phi(G_{i})$
there is a unique subchamber $(G_{i},R_{i})$ separating $P_{1}$ and $P_{2}$.
Note that since $\phi(G)$ is codimension $1$ in $\phi(F)$,
$\Aff(\phi(F)) \setminus S$ has two components $C_{1}$ and $C_{2}$, 
containing $P_{1}$ and $P_{2}$ respectively. Denote the 
map $\Lin(\phi(F)) \ra \Lin(\phi(F))/\Lin(S)$ by $\pi$.

Consider the weights of each $G_{i}$ which lie along the wall $F$, i.e.
\[
\alpha^{F}_{G_{i}} = \alpha_{G_{i}} \intersect (\Lin(\phi(F))/\Lin(S)).
\]
The weights $\alpha^{F}_{G_{i}} \subset \Lin(\phi(F))/\Lin(S) \iso \R$ 
divide into three classes: those which are zero, those lying
in $\pi(C_{1})$, and those lying in $\pi(C_{2})$.
Let 
\[
b_{i} = \# \{ \alpha_{G_{i},k} \mathrm{~lying~in~} \pi(C_{1}) \} 
\]
and
\[
f_{i} = \# \{ \alpha_{G_{i},k} \mathrm{~lying~in~} \pi(C_{2}) \}. 
\]
(We think of crossing the wall from $P_{1}$ to $P_{2}$ and $b_{i},f_{i}$ as
the number of weights pointing back and forward, respectively.) 
In practice, when we determine these numbers from a picture of the
X-ray, we can use the consistency condition. We
choose any vertex $v_{i}$ belonging to $G_{i}$ and count the
weights at $v$ which lie in $\phi(F)$ and which lie on the $C_{1}$ side
or the $C_{2}$ side respectively.

We will show that the signature (and the Poincar\'e polynomial and the
Euler characteristic) define invariants of the following form.
\begin{Def}\label{rec:invdef}
    A \textbf{recursive invariant} $I$ of a weighted X-ray $(\F,\phi,\alpha)$,
    with values in a ring $R$, consists of two pieces of data:
    \begin{enumerate}
      \item a map $I: Q_{\F} \ra R$, and
      \item a function $w_{I}: \N \times \N \ra R$, called 
            the \textit{wall-crossing function} of $I$.
    \end{enumerate}
    These must satisfy the following axioms for any wall $F \in \F$. 
    \begin{enumerate}
      \item \label{rec:insidecond} 
            Given two subchambers $(F,P_{1}),(F,P_{2})$
            separated by subchambers $(G_{i},R_{i})$, 
            let $b_{i},f_{i}$ be the number of weights of $G_{i}$ in $F_{i}$
            pointing toward $P_{1}$ and toward $P_{2}$ respectively. Then
            \begin{equation}\label{rec:invdefeq}
              I(P_{2}) - I(P_{1}) = \sum_{i} w_{I}(f_{i},b_{i}) I(R_{i}).
            \end{equation}
      \item \label{rec:outsidecond} 
            Given a subchamber $(F,P)$ adjacent to the boundary of $F$,
            let $(G,R)$ be the (necessarily unique) subchamber separating
            $P$ from the exterior of $F$. Let $f$ be the number of
            weights of $G$ in $F$ pointing into $F$ and $b$ be the number
            pointing out. Then
            \begin{equation}\label{rec:outsideeq}
              I(P) = w_{I}(f,b) I(R).
            \end{equation}
    \end{enumerate}
\end{Def}
Condition \ref{rec:outsidecond} is best understood as a special case
of condition \ref{rec:insidecond}, 
where one of the ``chambers'' is the exterior of the
wall $F$; in this view, it is understood that we always assign $I=0$ 
to this ``chamber.''\footnote{We avoid actually calling 
the exterior a chamber because it is neither
bounded nor convex; the trade-off is that we have to state 
condition \ref{rec:outsidecond} explicitly.}

The values of a recursive invariant on a given X-ray 
are completely determined by the function $w_{I}$ 
and the values of $I$ on the vertices. 
For, the values of $I$ on the subchambers of a given dimension $k$ 
are determined by its values on the subchambers of dimension $k-1$, 
along with $w_{I}$. This follows from the fact that we can get to any 
subchamber of dimension $k$ by starting outside the polytope and 
crossing a finite number of dimension $k-1$ walls. 

In fact, this gives an algorithm for calculating $I$
on all of the subchambers of $\F$, starting from the vertices and
working our way up. We give two examples 
in the case of the signature after Thm.\ \ref{rec:sigwallthm} below.

\subsection{Reduction to the Circle Case}\label{sec:circlereduction}

Note that any topological invariant $I$ of symplectic manifolds 
defines a function on the subchambers of a Hamiltonian X-ray by
\begin{equation}\label{rec:topinvonxray}
  I(F_{j},P) := I(M^{j}_{(q)}),\: \hbox{~for any~} q \in P.
\end{equation}
This invariant will be recursive in general if it is for circle
actions:
\begin{thm}\label{rec:circlereducthm}
  Let $I$ be a topological invariant of symplectic manifolds. Assume that
  the X-ray invariant defined by $I$ as in (\ref{rec:topinvonxray})
  is recursive 
  on the class of X-rays coming from Hamiltonian 
  circle actions. Then 
  this invariant is recursive on all Hamiltonian X-rays, with the
  same wall-crossing function.
\end{thm}

\begin{proof}
  We will show condition \ref{rec:insidecond}; condition
  \ref{rec:outsidecond} is very similar.
  Let $(M,\omega,\phi)$ be a Hamiltonian $T^{d}$-space. Let
  $(F,P_{1}),(F,P_{2})$ be adjacent subchambers of the X-ray of $M$.
  We note that everything takes place inside the wall $F$,
  which is an X-ray in its own right. Hence 
  we can assume that $M$ is connected and effective, $F=M$, 
  and $P_{1},P_{2}$ are chambers. 

  By Lemma \ref{rec:chamberoverlaplemma} the chambers $P_{1}, P_{2}$ 
  meet in a number of
  subchambers $(G_{j},R_{j})$ with common affine span $S$. Since the
  walls $(G_{j},\phi(G_{j}))$ are principal walls, the $G_{j}$ are 
  fixed point components of circles $T_{j} \subset T$.

  Let $a_{i} \in P_{i}$.
  Choose a $(d-1)$-torus $H \subset T$ such that under the projection
  $\pi: \ts \ra \hs$, $\pi(a_{1}) = \pi(a_{2})$. Let $L$ be the line 
  $\pi^{-1}(\pi(a_{i}))$. 
  Note that $L$ intersects the walls $\phi(G_{j})$ transversely.
  In  Fig.\ \ref{fig:stages} we illustrate the nonoverlapping case for
  simplicity; the projection $\pi$ is in the vertical 
  direction.
\begin{figure}[htbp]
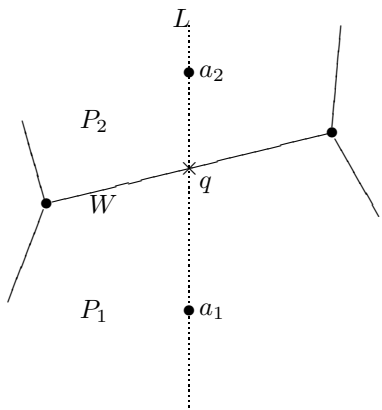

  \begin{center}
    \leavevmode
    \[ \xy/r3pc/:
  (0.5,0.125)="1", "1"+/r5ex/*{W},
  (3.5,0.875)="2", 
  (2,-2)="3", 
  (2,2)="4", "4"+/ul1ex/*{L},
  (2,1.5)="5", "5"+/r2ex/*{a_{2}},
  (2,-1)="6", "6"+/r2ex/*{a_{1}},
  (1,1)="7", "7"*{P_{2}},
  (1,-1)="8", "8"*{P_{1}},
  (2,0.5)="9", "9"*{\times}, "9"+/dr2ex/*{q},
  (0.1,-0.9)="10", 
  (0.25,1)="11", 
  (4,0)="12", 
  (3.6,2)="13", 
  "1"*{\bullet};                        
  "10" **\dir{-},                       
  "11" **\dir{-},                       
  "2" **\dir{-}, "2" *{\bullet};        
  "12" **\dir{-},                       
  "13" **\dir{-},                       
  "5"*{\bullet};                        
  "6"*{\bullet};                        
  "3";                                  
  "4" **\dir{.},                        
\endxy \]
    \caption{Reduction in Stages.}
    \label{fig:stages}
  \end{center}
\end{figure}

  Form the symplectic reduction $M_{H} = \phi^{-1}(L)$. 
  This has a residual action of the circle $T/H$, and the reduced
  spaces of $M_{H}$ by this circle action 
  at $a_{1}$ and $a_{2}$ are exactly the reduced spaces $M_{a_{1}}$ and
  $M_{a_{2}}$. Polarize $T/H$ so that the positive direction points
  from $a_{1}$ to $a_{2}$.

  The circle action on $M_{H}$ 
  has one fixed point component lying between $a_{1}$ and $a_{2}$
  for each wall $G_{j}$: these are exactly 
  $(G_{j})_{red} = M^{j}_{(q)}$, where $q = L \intersect \phi(G_{j})$.
  The weights of the circle action at these fixed point components
  are easily calculated.
  In particular, the positive weights of the circle action
  at $M^{j}_{(q)}$ 
  are in bijection with the forward weights of $G_{j}$, and the 
  negative weights are bijective with the backward weights. 

  Hence if the invariant $I$ obeys the wall-crossing formula 
  (\ref{rec:invdefeq}) for circle actions, then considering wall-crossing
  for $M_{H}$ gives
  \begin{equation}\label{cob:recinvstageseq}
    I(M_{a_{2}}) - I(M_{a_{1}}) = \sum_{j} w_{I}(f_{j},b_{j}) I(M^{j}_{(q)}).
  \end{equation}
  But by the definition of the recursive invariant associated to
  a topological invariant, this amounts to
  \begin{equation}\label{rec:invdefeqagain}
    I(P_{2}) - I(P_{1}) = \sum_{j} w_{I}(f_{j},b_{j}) I(R_{j}),
  \end{equation}
  which is the general wall-crossing formula.
\end{proof}

Hence any topological invariant which has a formula
in terms of weight data like that of the signature---and hence a wall-crossing
formula for circle actions---will be a recursive invariant.
We show below that this applies not only to the signature, but also to 
the Poincar\'e polynomial and the Euler characteristic.

\subsection{Application to the Signature}
\label{sec:recsig}

The signature is our first example of a recursive invariant.
Thm.\ \ref{rec:circlereducthm}  
combined with the circle wall-crossing formula (\ref{cob:wallcrossforcircle}) 
immediately gives
\begin{thm}\label{rec:sigwallthm}
  The signature $S(F_{j},P) = \Sign(M^{j}_{(q)})$, $q \in P$, 
  is a recursive invariant of Hamiltonian X-rays
  with wall-crossing function
\begin{equation}\label{sigwcfagain}
    w_{S}(f,b) = \left\{ \begin{array}{cr}
                       (-1)^{b} = -(-1)^{f} 
                               & \mathrm{~if~} f+b \mathrm{~is~odd} \\
                       0       & \mathrm{~if~} f+b \mathrm{~is~even}
                     \end{array}
             \right.     
  \end{equation}
\end{thm}
Note that the vertex function is just the signature
of the corresponding fixed point component,
\[
  S(p_{j}) = \Sign(F_{j}) \text{~for~} p_{j} = \phi(F_{j}), 
               F_{j} \subset M^{T}. 
\]
This and the wall-crossing function $w_{S}$ determine
the signature of all subreductions, and hence all regular reductions,
as discussed in the previous section.

\subsection{Simple Cases}\label{subsec:delzant}

In some cases, this recursive procedure collapses, making
the computation much simpler. First, we consider a trivial case.

Suppose $M$ is connected, the action of $T$ is effective,
and the dimension of $M$ is twice the rank of the torus.
Then the dimension of each regular reduced space is $0$,
and since a reduced space is connected, it is just a single
point. Hence the signature is always equal to $1$.

Such spaces are called \textit{symplectic toric varieties}; 
when they are smooth they are in fact the toric varieties of
complex algebraic geometry.
They are are characterized by the facts that they have isolated fixed points
and that their
X-rays have no internal structure, by the work of 
Delzant \cite{Delzant:1988HP}. That is, the only walls
of the X-ray are the faces of the moment polytope. 
Further, the moment polytopes of toric varieties are 
\textit{simple}, meaning that the edges coming out of each
vertex form a basis for $\ts$.

Now consider a general Hamiltonian $T^{d}$-space $M$
and look at its X-ray. Suppose that $M$ is connected
and the action is effective, for simplicity.
Suppose that the fixed points of $M$ are isolated.
Consider a stratum $F_{j}$ giving a wall
$\phi(F_{j})$ of dimension $r$. If the wall $\phi(F_{j})$
has no internal structure, we know that any
corresponding subreduction $M^{j}_{q}$ is just a point,
so $\Sign(M^{j}_{q})=1$. Hence the recursive procedure 
terminates early for this particular wall. We will
call such a wall, which corresponds to a toric 
stratum $F$, a \textit{Delzant wall} of the X-ray.

This can save a great deal of work in applying the recursive procedure
for calculating the signature of a reduction $M_{a}$. 

The simplest case is descibed by
\begin{lemma}\label{rec:whenalldelzantwalls}
  Let $M$ be a connected, effective Hamiltonian $T^{d}$-space
  with isolated fixed points,
  and let $(\F,\phi,\alpha)$ be its weighted X-ray.
  Then all principal walls of $\F$ are Delzant 
  exactly when for every vertex $p$, any $d$ weights at $p$ are independent. 
\end{lemma}
\begin{proof}
  First note that the condition that any $d$ weights are independent
  is equivalent to saying that any $(d-1)$ weights are independent
  and no other weight lies in their span.

  By the Darboux condition, a codimension $1$ wall is generated
  by a set of vectors which span a $(d-1)$-dimensional subspace
  of $\ts$. If the wall is Delzant, then it is a simple polytope, hence
  the weights in that subspace must be independent. So there must
  be only $(d-1)$ of them. Since no other weight can lie in their
  span, any $d$ weights will be independent.

  Conversely, if any $d$ weights are independent, then the walls 
  are generated by independent sets of $(d-1)$ vectors, hence 
  are Delzant.
\end{proof}

In such a case the recursive procedure terminates in
codimension $1$, giving a simple wall-crossing formula.

\begin{prop}\label{rec:wallcrossdelzantisolated}
  Let $M$ be a connected, effective Hamiltonian $T$-space
  with isolated fixed points. Assume that all of the
  principal walls are Delzant. Let $a_{1}$, $a_{2}$
  be regular values of the moment map, separated by
  walls $F_{1}$, \dots ,$F_{m}$, and associate 
  $b_{r},f_{r}$ to each wall as above. Then
\[
  \Sign(M_{a_{2}}) - \Sign(M_{a_{1}}) 
     = \sum_{r=1}^{m} w_{S}(f_{r},b_{r}). 
\]
\end{prop}

\begin{figure}[htbp]
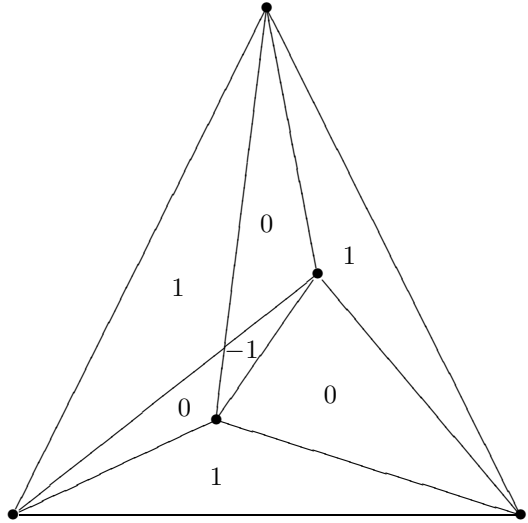

  \begin{center}
    \leavevmode
    \[ \xy/r4pc/:
  (0,0)="1",
  (4,0)="2", 
  (2,4)="3", 
  (1.6,0.75)="4", 
  (2.4,1.9)="5", 
  (1.3,1.8)="6", "6"*{1}, 
  (2,2.3)="7",     "7"*{0}, 
  (2.65,2.05)="8",   "8"*{1}, 
  (1.35,0.85)="9", "9"*{0}, 
  (1.8,1.3)="10", "10"*{-1}, 
  (2.5,0.95)="11", "11"*{0}, 
  (1.6,0.3)="12", "12"*{1}, 
  "1"*{\bullet};                        
  "2" **\dir{-}, "2" *{\bullet},        
  "3" **\dir{-}, "3" *{\bullet},        
  "4" **\dir{-}, "4" *{\bullet},        
  "5" **\dir{-}; "5" *{\bullet};        
  "2" **\dir{-},                        
  "3" **\dir{-},                        
  "4" **\dir{-};                        
  "2" **\dir{-},                        
  "3" **\dir{-};                        
  "2" **\dir{-};                        
\endxy \]
    \caption{Signatures of reduced spaces of $\CP^{4}$.}
    \label{fig:cobversioncp4}
  \end{center}
\end{figure}

\textit{Example.} Consider $\CP^{4}$ as a $T^{2}$-action, as in Sec.\ 
\ref{xraychap}. This space has isolated fixed points,
and all of the walls are Delzant (in particular, all multiplicities are
$1$). So the wall-crossing formula gives the signatures shown in Figure
\ref{fig:cobversioncp4}. For instance, we can calculate the
signature in the bottom chamber by crossing the bottom
wall from the outside. At either vertex of this wall, there are
three weights pointing up, none pointing down. Hence the difference
in signature is $W_{S}(3,0) = 1$. The other chambers are
calculated similarly.

This is not an isolated example. One can show (see \cite{Metzler:1997})
that if one starts with any symplectic toric variety and restricts
to a generic subtorus $T$, one obtains a Hamiltonian $T$-space
whose X-ray has only Delzant walls.

Now we turn to the non-Delzant case, where we do use the recursive
version of the formula.

\begin{figure}[htbp]
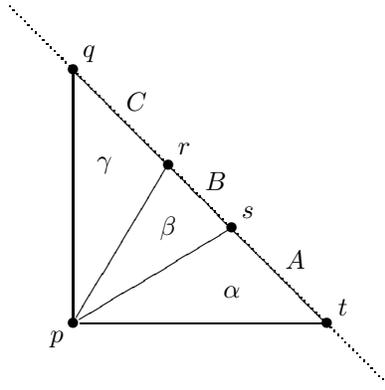

  \begin{center}
    \leavevmode
    \[ \xy/r2pc/:
  (0,0)="1", "1"+/dl2ex/*{p},
  (4,0)="2", "2"+/ur2ex/*{ t},
  (0,4)="3",  "3"+/ur2ex/*{ q},
  (1.5,2.5)="4", "4"+/ur2ex/*{ r},
  (2.5,1.5)="5", "5"+/ur2ex/*{ s},
  (-1,5)="6",
  (5,-1)="7",
  (2.5,0.5)="10", "10"*{\alpha},
  (1.5,1.5)="11", "11"*{\beta},
  (0.5,2.5)="12", "12"*{\gamma},
  (3.5,1)="15", "15"*{ A},
  (2.25,2.25)="16", "16"*{ B},
  (1,3.5)="17", "17"*{ C},
  "1"*{\bullet};                        
  "2" **\dir{-}, "2" *{\bullet},        
  "3" **\dir{-}, "3" *{\bullet},        
  "4" **\dir{-}, "4" *{\bullet},        
  "5" **\dir{-}; "5" *{\bullet};        
  "2";                                  
  "3" **\dir{-},                        
  "6";                                  
  "7" **\dir{.}, 
\endxy \]
    \caption{X-ray of $\CP^4$ with non-Delzant wall}
    \label{fig:nongenericcp4}
  \end{center}
\end{figure}

\textit{Example.\/} Consider $\CP^{4}$ as a $T^{2}$-space once again,
but this time with a non-generic projection, with resulting X-ray 
given in Figure \ref{fig:nongenericcp4}. One wall is now non-Delzant
(it is the projection of a $3$-simplex onto a line).
We could calculate the signatures in the chambers by crossing the
Delzant walls, but we will do it by crossing the non-Delzant
wall to show the idea of the recursive formula.

First we need the values on the vertices; but since all of the
$T$-fixed points are isolated, $S(p) = S(q) = \dots = 1$.

Next look at the diagonal wall. This is exactly the X-ray considered
in the example of Section \ref{sec:circle}. Hence the signatures
of the reductions are 
\[
  S(A) = S(C) = 1, \: S(B) = 0.
\]

Now we can cross the subchambers $A,B,C$ to find the signatures in the
chambers $\alpha,\beta,\gamma$. Crossing from the exterior of 
the polytope into $\alpha$
through $A$ gives 
\[
  S(\alpha) = w_{S}(0,1)\cdot S(A) = 1,
\]
and similarly $S(\gamma) = 1$. Crossing subchamber $B$ gives 
no change in signature, since $S(B)=0$, so $S(\beta) = 0.$

It is easy to see that the results on the chambers accord with what
we get by crossing the Delzant walls
emanating from the point $p$.


\subsection{The Poincar\'e Polynomial and the Euler Characteristic}
\label{sec:euler}

Theorem \ref{rec:circlereducthm} suggests that we look for
recursive formulas for other topological invariants. 
The cobordism approach which yielded Thm.\ \ref{rec:sigwallthm} will not
work for other invariants, because of Thm.\ \ref{cob:signatureisonly}.  
However, a little Morse theory 
leads to a wall-crossing formula for the Poincar\'e polynomial, 
and hence for the Euler characteristic. We will also note an 
intriguing connection to the signature formula. 

Thm.\ \ref{rec:circlereducthm} shows 
that any wall-crossing
result can be reduced to the case of a circle action. In that case, we can
derive a formula for the change in the Poincar\'e polynomial
as we cross a wall by the following procedure.

Let $(M,\omega,\phi)$ be a Hamiltonian $S^{1}$-space, with fixed
point components $F_{1},\dots F_{k}$. Identify $\mathrm{Lie}(S^{1})$
with $\R$ as usual. Pick one of the components
$F_{r}$, and let $\phi_{r} = \phi(F_{r})$. Assume for simplicity that
no other fixed point component is mapped under $\phi$ to $\phi_{r}$.
Let $a_{1},a_{2} \in \R$ be 
regular values of $\phi$ such that $a_{1} < \phi_{r} < a_{2}$, 
with no other $\phi(F_{j})$ lying between $a_{1}$ and $a_{2}$. 

We follow the suggesstion of Tolman (personal communication) 
in using the Morse theory of the square of two shifted moment maps, as
follows. 
Kirwan \cite{FK:1984CQ} showed 
that $\phi^{2}$ is an equivariantly perfect Morse
function.\footnote{Actually, $\phi^{2}$ is not necessarily Morse at
the critical value $0$, but $0$ is a minimum,
and hence the Morse theory still works. See \cite{FK:1984CQ} for a
detailed discussion.} 
Hence so are $(\phi-a_{1})^{2}$ and $(\phi-a_{2})^{2}$,
since the shifted maps are also $S^{1}$ moment maps. 
Note that the critical sets of $(\phi-a_{1})^{2}$ are the fixed points
of the action, plus the set $\phi^{-1}(a_{1})$; similarly for 
$(\phi-a_{2})^{2}$.

Calculating the equivariant Poincar\'e series $\tilde{P}_M$
of the manifold $M$ gives\footnote{We prefer the notation $P_{M}(t)$ for
the Poincar\'e polynomial but we often use P(M)(t) for easier reading when
there are many subscripts.}
\begin{eqnarray}\label{poincarecalc}
  \tilde{P}_{M}(t) &=& \sum_{F_{j}} P({F_{j}})(t) \frac{t^{\nu_{j}}}{1-t^{2}}
                   + \tilde{P}({\phi^{-1}(a_{1})})(t)\\
           &=& \sum_{F_{j}} P({F_{j}})(t) \frac{t^{\nu'_{j}}}{1-t^{2}}
                   + \tilde{P}({\phi^{-1}(a_{2})})(t)\\
\end{eqnarray}
where $\nu_{j}$ and $\nu'_{j}$ are the indices (the dimensions
of the negative normal bundles) of the fixed point set $F_{j}$
with respect to the Morse functions $(\phi-a_{1})^{2}$ and $(\phi-a_{2})^{2}$
respectively. Note that since $S^{1}$ acts freely on the level
sets $\phi^{-1}(a_{1})$ and $\phi^{-1}(a_{2})$, we have 
$\tilde{P}(\phi^{-1}(a_{1})) = P(M_{a_{1}})$ and similarly for $a_{2}$.

Since only the one fixed point set $F_{r}$ lies between the levels
$a_{1}$ and $a_{2}$, we have  
\[
  \nu_{j} = \nu'_{j} \hbox{~for~} j \ne r
\]
Furthermore, let $f$ be the number of positive weights in the
normal bundle to $F_{r}$ and $b$ be the number of negative weights.
Then $\nu_{r} = 2f$ and $\nu'_{r} = 2b$. (Each complex weight
gives two real dimensions.) When we equate the two right hand sides
in (\ref{poincarecalc}) above, all of the terms cancel except for
\begin{equation}\label{poincareresult}
   P({M_{a_{2}}})(t) - P({M_{a_{1}}})(t) = 
            P({F_{r}})(t) \frac{t^{2b}-t^{2f}}{1-t^{2}}.
\end{equation}

It is easy to see that if there are multiple fixed point
components $F_{1},\dots ,F_{m}$ 
with the same value of the 
moment map, we get
\begin{equation}\label{poincareresultoverlap}
   P({M_{a_{2}}})(t) - P({M_{a_{1}}})(t) = 
            \sum_{r=1}^{m} P({F_{r}})(t) 
            \frac{t^{2b_{r}}-t^{2f_{r}}}{1-t^{2}}.
\end{equation}

This is the wall-crossing formula for the circle case. 
Hence we immediately get a wall-crossing formula for a general
torus action, by Thm.~\ref{rec:circlereducthm}:
\begin{thm}\label{poincarewall}
  The X-ray invariant defined by the Poincar\'e polynomial, 
\[
P(F_{j},R) = P({M^{j}_{(q)}}) \in \Z[t], \quad  q \in R, 
\]
  is a recursive invariant of Hamiltonian X-rays
  with wall-crossing function 
  \begin{eqnarray}\label{rec:poincarewalleq}
     w_{P}(b,f) &=& \frac{t^{2b}-t^{2f}}{1-t^{2}}\\
                &=& t^{2f-2} + t^{2f-4}+ \dots + t^{2b} \:\: (f > b)\\
                &=& -t^{2b-2} - t^{2b-4}+ \dots - t^{2f} \:\: (b > f).
  \end{eqnarray}
\end{thm}

We get a wall-crossing formula for the Euler characteristic by setting
$t=-1$:
\begin{cor}\label{rec:eulerwall}
  The Euler characteristic $\chi$ defines a recursive invariant
  with wall-crossing function 
  \begin{equation}\label{rec:eulerwalleq}
    w_{\chi}(b,f) = b-f.
  \end{equation}
\end{cor}

\textit{Example.\/}
Consider again the X-ray of the nongeneric $2$-torus action on $\CP^{4}$
in Figure \ref{fig:nongenericcp4}. Once again we will start
with the non-Delzant wall to see how the recursive formula for
$P$ works. We evidently have $P(p) = P(q) = \dots = 1$. 

Looking at the diagonal wall, we cross $t$ to go from the outside of the
wall to subchamber $A$. The wall-crossing formula gives
\[
  P(A) = w_{P}(0,3) P(t) = 1 + t^{2} + t^{4}.
\]
(Note that this is $P(\CP^{2})$, as expected from the discussion on
p. \pageref{cob:examplecp3circle}.) Crossing from $A$ into $B$
through $s$ gives
\[
  P(B) = P(A) + w_{P}(1,2) P(s) = 1 + t^{2} + t^{4} + t^{2}.
\]
Crossing from $B$ into $C$ gives
\begin{align*}
  P(C) = P(B) +  w_{P}(2,1) P(r) &= 1 + 2 t^{2} + t^{4} - t^{2}\\
                                 &= 1 + t^{2} + t^{4}
\end{align*}
as we would expect from symmetry.

Now we cross the subchambers $A,B,C$ into the chambers 
$\alpha,\beta,\gamma$. Crossing from the exterior of 
the polytope into $\alpha$
through $A$ gives 
\[
  P(\alpha) = w_{P}(0,1) P(A) = 1 + t^{2} + t^{4},
\]
and similarly $P(\gamma) = 1 + t^{2} + t^{4}$. 
Crossing subchamber $B$ gives 
$P(\beta) =  1 + 2 t^{2} + t^{4}.$


Note that the wall-crossing functions (\ref{rec:eulerwalleq}) 
for the Euler characteristic and (\ref{sigwcfagain}) for 
the signature agree modulo $2$. This accords with the general
fact that $\Sign(M) \equiv \chi(M) \mod 2$ for all compact
oriented manifolds (\cite{HirBerJung:1992}).

It is an interesting fact (see \cite{Metzler:1997}) 
that the signature of a toric variety
is determined by its Poincar\'e polynomial---in fact 
  \[
    \Sign(V) = P_{V}(i) = P_{V}(\sqrt{-1})  
  \]
for $V$ a toric variety.
(This comes from the fact that the Hodge numbers $h^{p,q}$ vanish
for $p \ne q$.)

Comparing the wall-crossing formulas for the signature and the 
Poincar\'e polynomial shows that as far as the wall-crossing 
function is concerned, this still holds:
\[
  w_{S}(b,f) = w_{P}(b,f)(i).
\]

Since any recursive invariant is determined by its wall-crossing 
function and its values on vertices, we have the
\begin{thm}\label{eul:sigpoincare}
  Let $(M,\omega,\phi)$ be a Hamiltonian $T$-space. Assume that 
  the fixed point components $F \subset M^{T}$ are toric varieties. Then 
  \[
    \Sign(M^{j}_{(q)}) = P({M^{j}_{(q)}})(i) 
  \]
  for every $q \in \phi(F_{j})$, and in particular,
  \[
    \Sign(M_{a}) = P_{M_{a}}(i)
  \]
  for every $a \in \Delta_{reg}$.
\end{thm}
For example, if $M$ has isolated fixed points, the hypotheses of the theorem
are obviously satisfied.

Here is a simple consequence. Suppose that $N$ is a compact, connected
$4$-dimensional
orbifold that can be expressed as a reduced space for a Hamiltonian
torus action with isolated fixed points. Then 
\begin{align}\label{rec:4dim}
P_{N}(t) &= 1 + b_{1} t + b_{2} t^{2} + b_{1} t^{3} + t^{4} \\
\Sign(N) = P_{N}(i) &= 2 - b_{2}. 
\end{align}
If the diagonalized intersection form of $N$ has $p$ positive entries
and $n$ negative ones, we see that
\[
\Sign(N) = p - n = 2 - b_{2} 
\]
but $b_{2} = p + n$, so $p-n = 2 - p - n$ or
\[
p = 1.
\]
Hence any such $4$-manifold has exactly one 2-class with
positive self-intersection.





\end{document}